\documentclass{gtart}

\def\ifplaintex{\expandafter\ifx\csname documentclass\endcsname\relax}

\ifplaintex 
\else
\fi

\expandafter\ifx\csname beginpicture\endcsname\relax
\expandafter\ifx\csname documentclass\endcsname\relax
\input pictex \else
\input prepictex \input pictex \input postpictex \fi\fi

\def\gt{{\mathsurround=0pt\it $\cal G\mskip-2mu$eometry \&\ 
$\cal T\!\!$opology}}        

\def\gtp{{\mathsurround=0pt\it $\cal G\mskip-2mu$eometry \&\ 
$\cal T\!\!$opology $\cal P\!$ublications}}  

\def\lognumber#1{\def\thelognumber{#1}}
\def\volumenumber#1{\def\thevolumenumber{#1}}
\def\papernumber#1{\def\thepapernumber{#1}}
\def\volumeyear#1{\def\thevolumeyear{#1}}

\def\pagenumbers#1#2{\def\startpage{#1}\def\finishpage{#2}}
\def\published#1{\def\publishdate{#1}}
\def\proposed#1{\def\theproposer{#1}}
\def\seconded#1{\def\theseconders{#1}}
\def\received#1{\def\receiveddate{#1}}
\def\revised#1{\def\reviseddate{#1}}
\def\accepted#1{\def\accepteddate{#1}}
\def\asciititle#1{\def\theasciititle{#1}}

\def\coverauthors#1{\def\thecoverauthors{#1}}
\def\asciiauthors#1{\def\theasciiauthors{#1}}
\def\asciiaddress#1{\def\theasciiaddress{#1}}
\def\asciiemail#1{\def\theasciiemail{#1}}
\long\def\asciiabstract#1{\long\def\theasciiabstract{#1}}
\def\asciikeywords#1{\def\theasciikeywords{#1}}

\let\\\par\let\thelognumber\relax
\let\thevolumenumber\relax\let\thepapernumber\relax
\let\thevolumeyear\relax\let\thesamplenumber\relax\let\startpage\relax
\let\finishpage\relax\let\publishdate\relax\let\receiveddate\relax
\let\reviseddate\relax\let\accepteddate\relax\let\theasciititle\relax
\let\theasciiauthors\relax\let\theasciiaddress\relax
\let\theasciiabstract\relax\let\theasciikeywords\relax
\let\theasciiemail\relax\let\theshortauthors\relax\let\theshorttitle\relax
\let\thecoverauthors\relax

\long\def\maketitlep{   

\count0=\startpage

\gt\hfill      
\beginpicture
\setcoordinatesystem units <0.33truein, 0.33truein> point at 2.2 0.9
\setplotsymbol ({$\cal G$})
\plotsymbolspacing=9truept
\circulararc 315 degrees from 0 1 center at 0 0
\setplotsymbol ({$\cal T$})
\circulararc 315 degrees from 1 -1 center at 1 0
\endpicture
\break
{\small\ifx\thesamplenumber\relax 
Volume \else Sample
\fi\thevolumenumber\ (\thevolumeyear)
\startpage--\finishpage\nl
Published: \publishdate}
\vglue 0.5truein plus 0.4fil minus 0.1truein

{\parskip=0pt\leftskip 0pt plus 1fil\def\\{\par\smallskip}{\ifplaintex\large
\else\Large\fi\bf\thetitle}\par\medskip}   

\vglue 0pt plus 0.1fil 

{\parskip=0pt\leftskip 0pt plus 1fil\def\\{\par}{\sc\theauthors}
\par\medskip}

\vglue 0pt plus 0.1fil 

{\small\parskip=0pt\let\newline\\
{\leftskip 0pt plus 1fil\def\\{\par}{\sl\theaddress}\par}
\expandafter\ifx\theemail\relax    
\relax\else\vglue 5pt plus 0.02fil minus 2pt\def\\{\stdspace{\rm 
and}\stdspace} 
\cl{Email:\stdspace\tt\theemail}\fi
\ifx\theurl\relax                  
\relax\else\vglue 5pt plus 0.02fil minus 2pt\def\\{\stdspace{\rm 
and}\stdspace}
\cl{URL:\stdspace\tt\theurl}\fi\par}

\vglue 7pt plus 0.3fil minus 3pt

{\bf Abstract}
\vglue 5pt plus 0.1fil minus 2pt

\theabstract

\vglue 7pt plus 0.3fil minus 3pt

{\bf AMS Classification numbers}\quad Primary:\quad \theprimaryclass

Secondary:\quad \thesecondaryclass

\vglue 5pt plus 0.3fil minus 2pt

{\bf Keywords}\quad \thekeywords

\vglue 10pt plus 0.5fil minus 5pt

{\small  Proposed: \theproposer\hfill Received: \receiveddate\nl
Seconded: \theseconders\hfill 
\ifx\reviseddate\relax                         
Accepted: \accepteddate                        
\else
Revised: \reviseddate                          
\fi}
\eject
}       

\let\maketitlepage\maketitlep
\let\maketitle\maketitlepage

\font\phead=cmsl9 scaled 950
\font=cmsl9 scaled 1050
\font\pnum=cmbx10 scaled 913
\font\lnum=cmbx10 
\font\pfoot=cmsl9 scaled 950
\font=cmsl9 scaled 1050
\ifplaintex
\headline{\vbox to 0pt{\vskip -4.5mm\line{\small\phead\ifnum
\count0=\startpage ISSN 1364-0380 (on line)
1465-3060 (printed) \hfill {\pnum\folio}\else\ifodd\count0\def\\{ }%
\ifx\theshorttitle\relax\thetitle\else\theshorttitle\fi\hfill{\pnum\folio}
\else\def\\{ and }{\pnum\folio}\hfill\ifx\theshortauthors\relax\theauthors
\else\theshortauthors\fi\fi\fi}\vss}}
\footline{\vbox to 0pt{\vglue 0mm\line{\small\pfoot\ifnum\count0=\startpage
\copyright\ \gtp\hfill\else
\gt, Volume \thevolumenumber\ (\thevolumeyear)\hfill\fi}\vss
}}
\else
\makeatletter
\def\@oddhead{{\small\lhead\ifnum\count0=\startpage ISSN 1364-0380 (on line)
1465-3060 (printed) \hfill {\lnum\number\count0}\else\ifodd\count0
\def\\{ }\ifx\theshorttitle\relax \thetitle \else\theshorttitle\fi\hfill
{\lnum\number\count0}\else\def\\{ and }{\lnum\number\count0}
\hfill\ifx\theshortauthors\relax 
\theauthors\else\theshortauthors\fi\fi\fi}}\def\@evenhead{\@oddhead}
\def\@oddfoot{\small\lfoot\ifnum\count0=\startpage\copyright\ \gtp\hfill\else
\gt, Volume \thevolumenumber\ (\thevolumeyear)\hfill\fi}
\def\@evenfoot{\@oddfoot}
\makeatother
\fi

\newwrite\gtoutfile
\long\gdef\makeheadfile{  
{\def\\{, }\def\s{ }
\immediate\openout\gtoutfile head.xxx
\immediate\write\gtoutfile{Proxy-for: \ifx\theasciiauthors\relax
\theauthors\else\theasciiauthors\fi\s<\ifx\theasciiemail\relax\theemail\else\theasciiemail\fi>}
\immediate\write\gtoutfile{\noexpand\\}
\immediate\write\gtoutfile{Authors: \ifx\theasciiauthors\relax
\theauthors\else\theasciiauthors\fi}
{\def\\{ }\immediate\write\gtoutfile{Title: \ifx\theasciititle\relax
\thetitle\else\theasciititle\fi}}
\immediate\write\gtoutfile{Subj-class: GT or SG or MG etc}
\immediate\write\gtoutfile{MSC-class: \theprimaryclass\ifx\thesecondaryclass\relax\else, \thesecondaryclass\fi}
\immediate\write\gtoutfile{Journal-ref: Geom. Topol. \thevolumenumber
(\thevolumeyear) \startpage-\finishpage}
\immediate\write\gtoutfile{Comments: Published by Geometry and Topology at}
\immediate\write\gtoutfile{\s\s http://www.maths.warwick.ac.uk/gt/GTVol\thevolumenumber/paper\thepapernumber.abs.html}
\immediate\write\gtoutfile{\noexpand\\}
\immediate\write\gtoutfile{}
\ifx\theasciiabstract\relax
\immediate\write\gtoutfile{\theabstract}\else
\immediate\write\gtoutfile{\theasciiabstract}\fi
\immediate\write\gtoutfile{}
\immediate\write\gtoutfile{\noexpand\\}
\immediate\write\gtoutfile{}
\immediate\closeout\gtoutfile}}  

\def\maketitlepage{\maketitlep\makeheadfile}
\let\maketitle\maketitlepage

\lognumber{298}
\volumenumber{7}\papernumber{22}\volumeyear{2003}
\pagenumbers{773}{787}
\received{1 January 2003}
\revised{3 October 2003}
\published{13 November 2003}
\accepted{7 November 2003}
\proposed{Robion Kirby}
\seconded{Walter Neumann, Cameron Gordon}

\usepackage{rlepsf,amssymb,amsmath}

\def\S{Section }

\theoremstyle{plain}
\newtheorem*{theorem*}{Theorem}
\newtheorem*{proposition*}{Proposition}
\newtheorem*{corollary*}{Corollary}
\newtheorem{theorem}{Theorem}[section]
\newtheorem{lemma}[theorem]{Lemma}

\newtheorem{claim}[theorem]{Claim}

\theoremstyle{definition}
\newtheorem*{question*}{Question}
\newtheorem*{acknowledgements*}{Acknowledgements}

\newtheorem{example}[theorem]{Example}
\newtheorem{remark}[theorem]{Remark}
\numberwithin{equation}{section}
\numberwithin{figure}{section}

\font\nb=msbm10

\def \R{\hbox{\nb R}}

\def \Z{\hbox{\nb  Z}}

\def \Q{\hbox{\nb Q}}

\def \C{\hbox{\nb  C}}

\def \Hom{\hbox{Hom}}
\newcommand{\spinc}{{\rm{Spin}}$^c$}

\begin{document}

\title[Reidemeister--Turaev torsion modulo one]{Reidemeister--Turaev torsion 
modulo one\\of rational homology three--spheres}
\asciititle{Reidemeister-Turaev torsion 
modulo one of rational homology three-spheres}

\author{Florian Deloup\\Gw\'ena\"el Massuyeau}
\coverauthors{Florian Deloup\\Gw\noexpand\'ena\noexpand\"el Massuyeau}
\asciiauthors{Florian Deloup\\Gwenael Massuyeau}

\address{Laboratoire Emile Picard, UMR 5580 CNRS/Univ. Paul Sabatier \\
118 route de Narbonne, 31062 Toulouse Cedex 04, France}
\secondaddress{Laboratoire Jean Leray, UMR 6629 CNRS/Univ. de Nantes \\
2 rue de la Houssi\-ni\`e\-re, BP 92208, 44322 Nantes Cedex 03, France}

\email{deloup@picard.ups-tlse.fr}
\secondemail{massuyea@math.univ-nantes.fr}

\asciiemail{deloup@picard.ups-tlse.fr, massuyea@math.univ-nantes.fr}

\asciiaddress{Laboratoire Emile Picard, UMR 5580 CNRS/Univ. Paul
Sabatier\\118 route de Narbonne, 31062 Toulouse Cedex 04, 
France\\and\\Laboratoire Jean Leray, UMR 6629 CNRS/Univ. de 
Nantes\\2 rue de la Houssiniere, BP 92208, 44322 Nantes Cedex 03, France}

\primaryclass{57M27} 
\secondaryclass{57Q10, 57R15}

\keywords{Rational homology $3$--sphere, Reidemeister torsion, complex spin structure, quadratic function}
\asciikeywords{Rational homology 3-sphere, Reidemeister torsion, complex spin structure, quadratic function}

\begin{abstract} 
Given an oriented rational homology $3$--sphere $M$, it is known how to associate to any Spin$^c$--structure $\sigma$
on $M$ two quadratic functions over the linking pairing. 
One quadratic function is derived from the reduction modulo $1$ of the Reidemeister--Turaev torsion of 
$(M,\sigma)$, while the other one can be defined using the intersection pairing 
of an appropriate compact oriented $4$--manifold with boundary $M$.

In this paper, using surgery presentations of the manifold $M$,
we prove that those two quadratic functions coincide. Our proof relies 
on the comparison between two distinct combinatorial descriptions 
of Spin$^c$--structures on $M$: Turaev's charges vs Chern vectors.
\end{abstract}

\asciiabstract{Given an oriented rational homology 3-sphere M, it is
known how to associate to any Spin^c-structure \sigma on M two
quadratic functions over the linking pairing.  One quadratic function
is derived from the reduction modulo 1 of the Reidemeister-Turaev
torsion of (M,\sigma), while the other one can be defined using the
intersection pairing of an appropriate compact oriented 4-manifold
with boundary M.  In this paper, using surgery presentations of the
manifold M, we prove that those two quadratic functions
coincide. Our proof relies on the comparison between two distinct
combinatorial descriptions of Spin^c-structures on M Turaev's
charges vs Chern vectors.}

\maketitlepage

\section{Introduction and statement of the result}
\subsection{Introduction}

Any closed oriented $3$--manifold $M$ can be equipped with a \emph{complex
spin structure}, or \emph{Spin$^c$--structure}. While they seem to have been
originally introduced in the '50s and '60s \cite{Hi}, in the framework
of Dirac operators and $K$--theory \cite{LM}, the revival of interest in 
Spin$^c$--structures over the last decade is certainly due to symplectic
geometry and Seiberg--Witten invariants of $4$--manifolds. For a general
introduction to Spin$^c$--structures, the reader is referred to \cite{LM}. 
It was observed somewhat more recently \cite{TSpinc} that, 
in dimension 3, Spin$^c$--structures have a simple and natural interpretation: 
any Spin$^c$--structure on a closed oriented $3$--manifold $M$ 
can be represented by a nowhere vanishing vector field on $M$. 
This enabled Turaev to reinterpret a topological invariant of
Euler structures on $3$--manifolds, which he had introduced earlier,
 as an invariant of Spin$^c$--structures. Since this invariant
is a refinement of the Reidemeister torsion, we call this invariant
the Reidemeister--Turaev torsion.

We will be interested in the restriction of this invariant to the class of
rational homology $3$--spheres. Our work is motivated by and based on 
two observations.
\begin{enumerate}
\item[-] On the one hand, there is the following special feature of the
Reidemeister--Turaev torsion
 $\tau_{M, \sigma}$ of an oriented rational homology $3$--sphere $M$ with a
Spin$^c$--structure $\sigma$: its reduction modulo 1
induces a quadratic function $q_{M,\sigma}$ over the linking pairing 
$\lambda_{M}$ \cite{TTorsions}. 

\item[-] On the other hand, there is a canonical bijective
correspondence, denoted by $\sigma \mapsto \phi_{M,\sigma}$, between
\spinc--structures on $M$ and quadratic functions over the linking
pairing $\lambda_{M}$ \cite{LW, Gi, De}. The quadratic function $\phi_{M,\sigma}$
can be defined, extrinsically, using the intersection pairing of a compact
oriented $4$--manifold with boundary $M$ and first Betti number equal to zero.
\end{enumerate}
Thus, the question naturally arises to compare the quadratic functions
$q_{M,\sigma}$ and $\phi_{M,\sigma}$. 

\subsection{Statement of the result}

Let us begin by developing the above two observations and fixing some 
notations.

The Reidemeister--Turaev torsion of a closed oriented
$3$--manifold equip\-ped with a Spin$^c$--structure is a fundamental
topological invariant. A concise and almost self-contained introduction
is \cite{TBDZ}. A broader introduction is \cite{TBook}, 
while the monographs \cite{Nic,TTorsions} contain the
most recent developments. We give here a succinct presentation sufficient
for our purpose. 

Let $M$ be a connected oriented $3$--manifold, compact without
boundary. All homology and cohomology groups will be with integral
coefficients unless explicity stated otherwise. We set $H =
H_{1}(M)$, the first homology group, written multiplicatively. Let
$Q(H)$ denote the classical ring of fractions of the group ring $\Z[H]$.
The \emph{maximal Abelian Reidemeister
 torsion} $\tau(M)$ of $M$ is an element in $Q(H)$ defined up to
multiplication by an element of  $\pm H\subset Q(H)$. This invariant, 
defined in \cite{TDebut}, can be thought of as a generalization of the Alexander polynomial.
Next, its indeterminacy in $\pm H$ can be disposed of by
specifying two extra structures: a homology orientation of $M$ and
an Euler structure of $M$ (see \cite{Tu89}).
On the one hand, using the intersection pairing, the choosen orientation of $M$
induces a canonical homology orientation. On the other hand, the Euler 
structures on $M$, defined as punctured homotopy classes of nowhere vanishing vector fields on $M$, 
are in canonical bijective correspondence with the \spinc--structures on $M$ \cite{TSpinc}. Therefore,
if $(M,\sigma)$ is a connected closed \spinc--manifold of dimension $3$, 
one can define its \emph{Reidemeister--Turaev torsion}
$$\tau(M,\sigma)\in Q(H).$$ 
It has the following equivariance property:
\begin{equation}
\label{eq:equivariance_of_torsion}
\forall h\in H,\ \ \ h\cdot\tau(M,\sigma)=\tau(M,h\cdot \sigma)\in Q(H).
\end{equation}
Here, the left hand side involves a multiplication in $Q(H)$ while, in the
right hand side, $h \cdot \sigma$ involves the free and
transitive action of $H^2(M)$ (or $H_{1}(M)$ via Poincar\'e duality)
on the set Spin$^{c}(M)$: see, eg, \cite{LM}.

Now and throughout the paper, {\bf{we assume that $M$ is an
oriented rational homology $3$--sphere}}, ie, we suppose that 
$$H_*(M;\Q)=H_*\left(\mathbf{S}^3;\Q\right).$$ 
Then $H$ is finite and $Q(H)=\Q[H]$. Hence $\tau(M,\sigma)$ determines a 
function
$\tau_\sigma: H\to \Q$ such that 
\begin{displaymath}
\tau(M,\sigma)=\sum_{h\in H} \tau_\sigma(h)\cdot h \in \Q[H].
\end{displaymath}
It has been proved in \cite[Theorem 4.3.1]{TSpinc} that the modulo $1$
reduction of the function $\tau_{\sigma}$ satisfies the property that
\begin{equation}
\label{eq:linking_and_torsion}
\forall h_1, h_2 \in H, \quad
\tau_{\sigma}(h_1 h_2) - \tau_{\sigma}(h_1)
- \tau_{\sigma}(h_2) + \tau_{\sigma}(1)
= - \lambda_{M}(h_1, h_2)\ \ \hbox{mod}\ 1.
\end{equation}
Here, $\lambda_M: H \times H \to \Q/\Z$ denotes the \emph{linking pairing} 
of $M$:
this is a symmetric nondegenerate bilinear pairing, which gives partial 
information on the way
knots are linked in the manifold $M$ \cite{ST}. It immediately follows from (\ref{eq:linking_and_torsion}) that
$$\forall h\in H, \quad \tau_{\sigma}(h) = \tau_{\sigma}(1) - q_{M,\sigma}\left(h^{-1}\right) \ \hbox{mod}\ 1,$$
where $q_{M,\sigma}$ is a \emph{quadratic function over} the linking pairing
$\lambda_{M}$, in the sense that it satisfies the following property:
$$
\forall h, k \in H, \ \ q_{M,\sigma}(hk) - q_{M,\sigma}(h) - q_{M,\sigma}(k) = \lambda_M(h,k).
$$  
It is also easily seen from (\ref{eq:equivariance_of_torsion}) and (\ref{eq:linking_and_torsion}) that
\begin{equation}
\label{eq:qM_is_affine} \forall h \in H, \quad q_{M,h\cdot
\sigma}=q_{M,\sigma}+ \lambda_M\left(h,-\right).
\end{equation}
This equation suggests to define the following free transitive action of the 
group $H$ on the set Quad$(\lambda_{M})$ of quadratic functions over $\lambda_{M}$:
$$
H \times {\rm{Quad}}(\lambda_{M}) \to {\rm{Quad}}(\lambda_{M}), \ (h, q)  \mapsto  h \cdot q
$$ 
where 
$$
\forall x \in H, \ (h \cdot q)(x) = q(x) +\lambda_{M}\left(h, x\right).
$$

On the other hand, it is known \cite{LW,Gi,De} (see \cite{DM} for 
arbitrary closed oriented $3$--manifolds) how to define another bijective $H$--equivariant correspondence
$$
\hbox{Spin}^{c}(M) \to \hbox{Quad}(\lambda_{M}), \
\sigma \mapsto \phi_{M,\sigma}.
$$
This map is defined combinatorially, starting from a surgery presentation of the manifold $M$ 
and using its linking matrix. (The detailed construction will be recalled 
in subsection \ref{subsec:quad}.)\\

\begin{theorem*}\
For any oriented rational homology $3$--sphere $M$ equipped with a Spin$^c$--structure $\sigma$,
the quadratic functions $q_{M,\sigma}$ and $\phi_{M,\sigma}$ are equal.
\end{theorem*}

In his monograph \cite{Nic}, Nicolaescu has proved
the same result, with an analytic proof based on the connection
between the Reidemeister--Turaev torsion and the Seiberg--Witten invariant.
Our proof is combinatorial and purely topological. A surgery presentation of $M$ 
provides two combinatorial descriptions of Spin$^c$--structures on $M$. 
One description (called \emph{charges}) is defined by Turaev in \cite{TSurg} 
in terms of the complement in ${\mathbf{S}}^{3}$ 
of the framed surgery link, and is used there to compute $\tau(M,\sigma)$. 
Another description (called \emph{Chern vectors}) relies on the $4$--manifold with boundary $M$ 
associated to the surgery presentation, and is well suited for the computation of $\phi_{M,\sigma}$. 
Our main contribution consists in comparing those two descriptions of Spin$^c$--structures.\\

Before going into the proof of the Theorem,
let us discuss the following immediate consequence.

\begin{corollary*}\
The quadratic function $\phi_{M,\sigma}$ is determined by $\tau(M,\sigma) \ {\rm{mod}}\ 1$.
\end{corollary*}

We claim that the converse of the Corollary does not hold. To justify this, define the  ``constant''
$$
c_\sigma=\tau_\sigma(1) \ \textrm{mod} \ 1.
$$
From (\ref{eq:equivariance_of_torsion}), we obtain that
\begin{equation}
\label{eq:c} 
\forall h\in H,\ \ \ c_{h\cdot\sigma}=c_\sigma-\phi_{M,\sigma}\left(h\right).
\end{equation}
Let also $d_\sigma\in \R/\Z$ be such that
$$
\exp\left(2i\pi\ d_\sigma\right)= \frac{1}{\sqrt{|H|}}\cdot
\sum_{x\in H} \exp\left(2i\pi\ \phi_{M,\sigma}(x)\right)\in \C.
$$
Since $\phi_{M,\sigma}$ is nondegenerate, the Gauss sum on the right
hand side is well-known to be a complex number of modulus $1$. It can also be proved
that $d_{\sigma} \in \Q/\Z$. Observe that
\begin{equation}
\label{eq:d}
d_{h\cdot\sigma}=d_\sigma-\phi_{M,\sigma}\left(h\right).
\end{equation}
As an immediate consequence of (\ref{eq:c}) and (\ref{eq:d}), we obtain the following

\begin{proposition*}\
The number $c(M)=c_\sigma-d_{\sigma}\in \Q/\Z$
is a topological invariant of the oriented rational homology $3$--sphere $M$.
\end{proposition*} 

Explicit computations can be performed on the lens spaces. For
instance, we find that $8c\left(L(7,1)\right)=3/7\neq 2/7= 8
c\left(L(7,2)\right)$; since $L(7,1)$ and $L(7,2)$ have isomorphic
linking pairings, we deduce that $c(M)$ can not be computed from 
$\phi_{M,\sigma}$.

It is not difficult to verify that $c(M)$ is additive under connected sums, 
vanishes if $M$ is an integer homology $3$--sphere and changes sign when the
orientation of $M$ is reversed. Let $\lambda(M)\in \Q$ denote the
Casson-Walker invariant of $M$
in Lescop's normalization \cite{Les}. We ask the following

\begin{question*}\
Does the invariant $c(M)\in \Q/\Z$ coincide with $-\lambda(M)/|H| \ \textrm{mod} \ 1$?
\end{question*}

\begin{acknowledgements*}\
The first author is an EU Marie Curie Research Fellow (HPMF 2001--01174) at the 
Einstein Institute of Mathematics, the Hebrew University of Jerusalem.
\end{acknowledgements*}

\section{Chern vectors and charges}

\label{sec:preliminaries}

This section contains preliminary material for the proof of the Theorem
(\S \ref{sec:proof}). The heart of this section is devoted to the presentation
of two equivalent, but distinct, combinatorial descriptions of
complex spin structures on $M$. The
proof of this equivalence will be given in  \S \ref{sec:proof}.
Even though we shall not need it, note that subsections
\ref{subsec:surgery}, \ref{subsec:Chern-vectors} and
\ref{subsec:charges} are valid for \emph{any} closed oriented connected
$3$--manifold (ie, with arbitrary first Betti number).

As a convention, boundaries of oriented manifolds will be always given orientation
by the ``outward normal vector first'' rule.

\subsection{Surgery presentation} \label{subsec:surgery}

In this paragraph and throughout \S \ref{sec:preliminaries},
we fix an ordered oriented framed $n$--component link $L$
in $\mathbf{S}^{3}$, such that the oriented $3$--manifold $V_L$ obtained
from $\mathbf{S}^{3}$ by surgery along $L$ is diffeomorphic 
to our oriented rational homology $3$--sphere $M$.

\begin{figure}
\centerline{\relabelbox \small
\epsfxsize 4truein \epsfbox{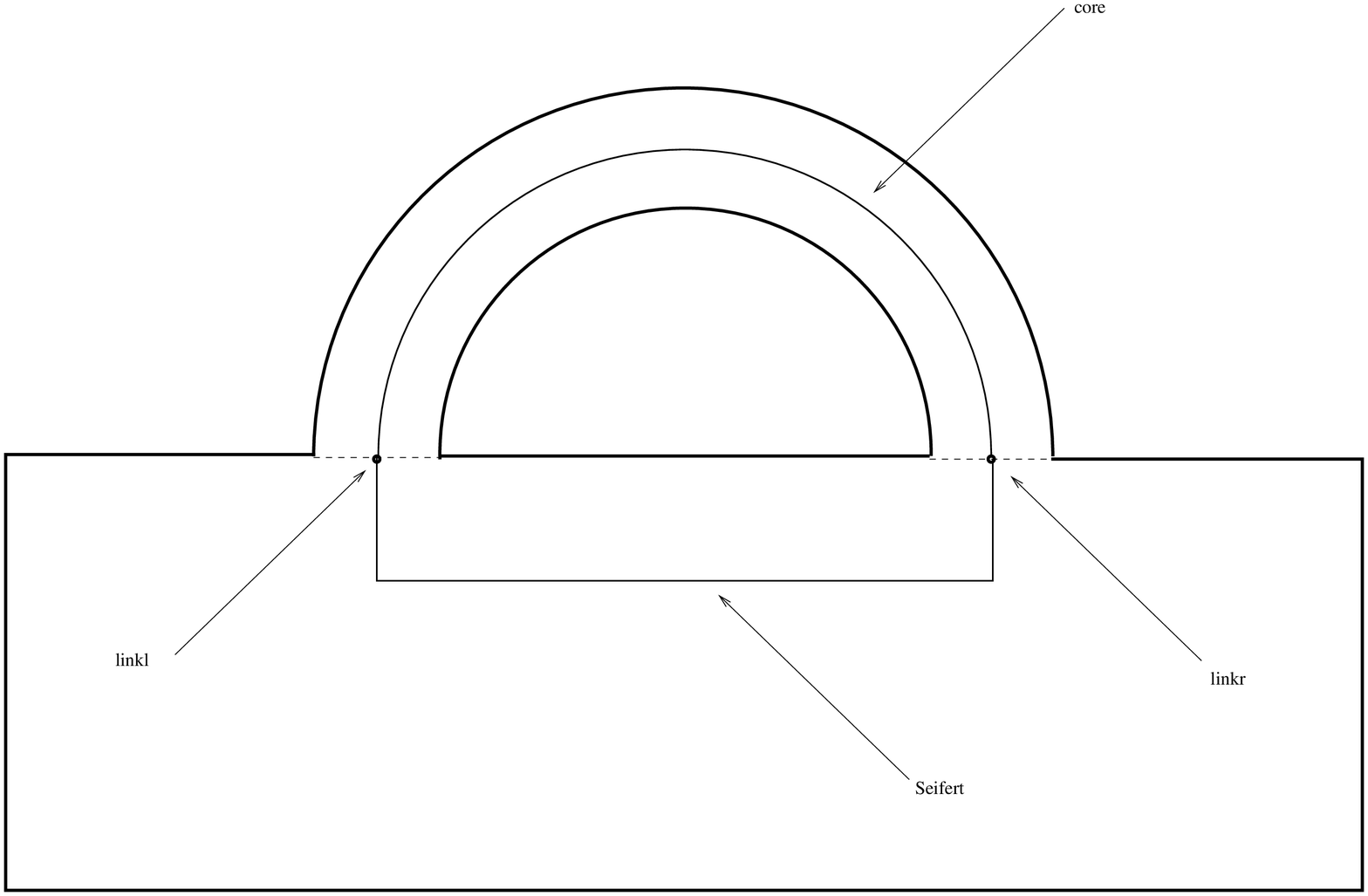}
\relabel {Seifert}{$-\Sigma_i$}
\relabel {core}{$\left(\mathbf{D}^2\times 0\right)_i$}
\relabel {linkl}{$L_i$}
\relabel {linkr}{$L_i$}
\endrelabelbox}
\caption{The preferred basis of $H_2\left(W_L\right)$}
\label{fig:handle}
\end{figure}

Let $b_{ij} = {\rm{lk}}_{\mathbf{S}^{3}}(L_{i},L_{j})$
for all $1 \leq i \not= j \leq n$, and let $b_{ii}$ be the framing number of
$L_{i}$ for all $ 1 \leq i \leq n$. We denote by $B_{L}=
(b_{ij})_{i,j=1,\dots,n}$ the linking matrix of $L$ in ${\mathbf{S}}^{3}$.
We also denote by $W_{L}$ the
\emph{trace} of the surgery. In other words,
$$ M = V_L=\partial W_L \quad \hbox{with} \quad
W_L=\mathbf{D}^4\cup \bigcup_{i=1}^n
\left(\mathbf{D}^2\times \mathbf{D}^2\right)_i, $$
where the $2$--handle $\left(\mathbf{D}^2\times \mathbf{D}^2\right)_i$
is attached by embedding $-\left(\mathbf{S}^1\times
\mathbf{D}^2\right)_i$ into $\mathbf{S}^3=\partial \mathbf{D}^4$
in accordance with the specified framing and orientation of $L_i$.
The group $H_{2}(W_{L})$ is free Abelian of rank $n$. It is given
the \emph{preferred} basis $([S_{1}],\dots,[S_{n}])$.
Here, the closed surface $S_i$ is taken to be
$$
S_i=\left(\mathbf{D}^2\times 0\right)_i
\cup \left(-\Sigma_i\right), 
$$
where $\Sigma_i$ is a Seifert surface
for $L_i$ in $\mathbf{S}^3$ which has been pushed into
the interior of $\mathbf{D}^4$ as shown in Figure \ref{fig:handle}.
Also, $H^2(W_L)$ will be identified with
${\Hom}(H_2(W_L),\Z)$ by Kronecker evaluation,
and will be given the dual basis. Note that the matrix of 
the intersection pairing $\bullet: H_{2}(W_L) \times H_2(W_L) \to \Z$
relatively to the preferred basis of $H_{2}(W_L)$ is $B_L$.

\subsection{Chern vectors} \label{subsec:Chern-vectors}

We define the set of {\it{Chern vectors}} (associated to the link $L$) to be
\begin{displaymath}
{\tilde{\mathcal{V}}}_L=
\left\{ s=(s_i)_{i=1}^n\in\Z^n \ : \
\forall i=1,\dots,n,\ s_i \equiv b_{ii}\textrm{ mod }2\right\}.
\end{displaymath}
Set $\displaystyle {\mathcal{V}}_L =
\frac{{\tilde{\mathcal{V}}}_L}{ 2\cdot {\rm{Im}}\ B_L}$.
A basic result of \cite{DM} (where the reader is referred to for
full details) asserts that
\begin{equation}
 {\rm{Spin}}^{c}(V_{L}) \simeq {\mathcal{V}}_L. \label{eq:DM-iso}
\end{equation} This is our first combinatorial description of
Spin$^c$--structures on $V_{L}$, which we now recall briefly. Let $\sigma \in {\rm{Spin}}^{c}(V_{L})$.
Extend $\sigma$ to a Spin$^c$--structure $\tilde{\sigma}
\in {\rm{Spin}}^{c}(W_{L})$. Thus the Chern class
$c(\tilde{\sigma}) \in H^{2}(W_{L}) \simeq$ Hom$(H_{2}(W_L), \Z)$ is given by
an element in $\Z^{n}$
(according to the basis dual to the preferred basis).
The isomorphism (\ref{eq:DM-iso}) is induced by the
map $\sigma \mapsto c({\tilde{\sigma}})$.

\subsection{Charges} \label{subsec:charges}

Charges were introduced by Turaev in \cite{TSurg}, 
as a combinatorial description of Euler structures. We give a brief description.

The set of
{\it{charges}} (associated to the link $L$) is defined to be

\begin{displaymath}
{\tilde{\mathcal{C}}}_L=
 \left\{ k=(k_i)_{i=1}^n\in\Z^n \ : \
\forall i=1,\dots,n,\ k_i \equiv 1 + \sum_{{1 \leq j \leq n},\ 
{j \not= i}} b_{ij} \textrm{ mod }2 \right\}.
\end{displaymath}
Set $\displaystyle {\mathcal{C}}_L= \frac{{\tilde{\mathcal{C}}}_L}{2\cdot
{\rm{Im}}\ B_L}$. We shall recall below that
\begin{equation}
  {\rm{Spin}}^{c}(V_{L})  \simeq  {\mathcal{C}}_L. \label{eq:T-iso}
\end{equation}

We can alternatively view $V_{L}$, without reference to $W_{L}$, as
$$ V_{L} = {\mathbf{E}} \cup \bigcup_{i = 1}^{n} Z_{i},$$
where ${\mathbf{E}}$ denotes the exterior of a tubular neighborhood of
$L$ in ${\mathbf{S}}^{3}$
and $Z_{i}$ is a (reglued) solid torus,
homeomorphic to ${\mathbf{S}}^{1} \times {\mathbf{D}}^{2}$. A solid
torus $Z$ is said to be \emph{directed}
when its core is oriented. We direct the solid torus $Z_j$ in the
following way: we denote by $m_j\subset {\mathbf{E}}$ the meridian of $L_j$
which is oriented so that  $\hbox{lk}_{\mathbf{S}^{3}}(m_j, L_j) = +1$,
and we require the oriented core of $Z_j$ to be isotopic in $V_L$ to $m_j$.

In general, let $N$ be a compact oriented $3$--manifold with boundary $\partial
N$ endowed,  this time, with a Spin--structure $\sigma$. There is a well-defined
set of \emph{Spin$^c$--structures on $N$ relative to $\sigma$}, 
denoted by $\hbox{Spin}^c(N,\sigma)$.
The Abelian group $H^2(N,\partial N)$ acts freely and  transitively on
$\hbox{Spin}^c(N,\sigma)$. Also, there is a \emph{Chern class map}
$$c: \hbox{Spin}^c(N,\sigma) \to H^2(N,\partial N)$$
which is affine over the square map (where $H^{2}(N,\partial N)$ is
written multiplicatively). For details about relative
\spinc--structures and their gluings, see \cite{DM}.

The torus ${\mathbf{S}}^{1} \times {\mathbf{S}}^{1}$ has
a canonical Spin--structure $\sigma^{0}$, which is induced by its Lie
group structure. Hence $\partial {\mathbf{E}}$ can be endowed with a
distinguished Spin--structure, which is denoted by $\cup_{i=1}^{n} \sigma^0$.
A directed solid torus $Z$ has a \emph{distinguished}
\spinc--structure relative to the canonical
Spin--structure $\sigma^0$ on $\partial Z$:
this is the one whose Chern class
is Poincar\'e dual to the opposite of the oriented core of $Z$.
Hence by gluing any Spin$^c$--structure on ${\mathbf{E}}$
relative to $\cup_{i=1}^{n} \sigma^0$
to the distinguished relative \spinc--structures
on the directed solid tori $Z_j$'s, we define a map
$$g: \hbox{Spin}^c\left({\mathbf{E}}, {\cup}_{i=1}^{n} \sigma^0\right) \to
{\hbox{Spin}}^c(V_{L}).$$
This map $g$ is affine, via the Poincar\'e duality isomorphisms
$P: H_{1}({\mathbf{E}}) \to H^{2}(\mathbf{E}, \partial {\mathbf{E}})$ and
$P: H_{1}(V_{L}) \to H^{2}(V_{L})$, over the natural inclusion homomorphism
$H_1({\mathbf{E}})\to H_1(V_L)$. In particular, $g$ is onto.

Another useful general fact is that the Chern class $c(\alpha)$ of a Spin$^c$--structure
$\alpha$ relative to a Spin--structure on the boundary has a nice explicit
expression modulo 2, which we briefly explain. Let $S$ be a closed oriented
surface. Denote by Quad$(S)$ the set of quadratic functions over the mod $2$ intersection pairing of $S$.
Hence, an element $q \in $ Quad$(S)$ is a map $q:H_{1}(S;\Z_{2}) \to \Z_{2}$ such that
$q(x + y) - q(x) - q(y) = x\bullet y$ for all $x, y \in H_{1}(S;\Z_{2})$, where $\bullet$
denotes the mod $2$ intersection pairing. 
The Atiyah-Johnson correspondence \cite{A,J} is a bijective $H_{1}(S,\Z_{2})$--equivariant map 
$$J: {\rm{Spin}}(S) \to {\rm{Quad}}(S), \ \sigma \mapsto J_{\sigma}.$$
Here, the function $J_\sigma$ is defined, for any simple oriented closed curve $\gamma$,
by $J_{\sigma}([\gamma]) = 1$ or $0$ according to whether
$(\gamma, \sigma|_{\gamma})$ is homotopic to ${\mathbf{S}}^{1}$
with the Spin--structure induced from the Lie group structure or not \cite[pages 35--36]{K}.

\begin{lemma}[See \cite{DM}]
\label{lem:Chern_mod2}
Let $N$ be a compact oriented $3$--manifold with boundary,
$\sigma \in {\rm{Spin}}(\partial N)$
and $\alpha\in {\rm{Spin}}^c(N,\sigma)$. Then
\begin{displaymath}
\forall y\in H_2(N,\partial N), \quad
\langle c(\alpha),y\rangle \equiv J_\sigma\left(\partial_*(y)\right)
{\rm{mod}}\ 2,
\end{displaymath}
where $\langle\cdot,\cdot\rangle$ denotes Kronecker evaluation,
and where $\partial_*:H_2(N,\partial N) \to H_1(\partial N)$
is the connecting homomorphism of the pair $(N,\partial N)$.
\end{lemma}

A canonical bijection between $\hbox{Spin}^c\left({\mathbf{E}},{\cup}_{i=1}^{n} \sigma^0\right)$
and ${\tilde{{\mathcal{C}}}}_{L}$ can be defined in the following way: for any
$\alpha \in \hbox{Spin}^c\left({\mathbf{E}}, {\cup}_{i=1}^{n} \sigma^0\right)$,
calculate $P^{-1}c(\alpha) \in H_1({\mathbf{E}})$ and identify
$H_1({\mathbf{E}})$ with $\Z^n$
taking the meridians $\left([m_{1}], \ldots, [m_n]\right)$ as a basis;
it is a consequence of Lemma \ref{lem:Chern_mod2}
that the multi-integer we obtain
is actually a charge on $L$. Thus, since $g$ is surjective and
since $\hbox{Ker}\left(H_1({\mathbf{E}})\to H_1(V_L)\right)$ is generated
by the $n$
characteristic curves of the surgery,
it follows that the map $g$ induces a bijection
$$\frac{{\tilde{\mathcal{C}}}_{L}}{2\cdot\hbox{Im}\ B_{L}}\to
\hbox{Spin}^c(V_{L})$$
as claimed.\\

\subsection{The quadratic function $\phi_{M,\sigma}$} \label{subsec:quad}

In this paragraph, we recall how to compute the quadratic function
$\phi_{M,\sigma}$ \cite{LW,Gi,De} from the surgery presentation $L$ for $M$ 
and a Chern vector $s \in\Z^{n}$ representing $\sigma \in {\rm{Spin}}^{c}(M)$. By the
homology exact sequence associated to the pair $(W_L, V_L)$, the
choice of the preferred basis for $H_2(W_L)$ induces an
identification
\begin{equation}
H \simeq {\rm{Coker}}\ B_{L} = \Z^n / {\hbox{Im}}\ B_{L}.
\label{eq:identification}
\end{equation}
Let $x \in H$ and let $X \in \Z^{n}$ be a representative of $x$ 
by $(\ref{eq:identification})$. We have
\begin{equation}
\phi_{M,\sigma}(x) =  -\frac{1}{2}
\left(X^{\rm{T}}\cdot {B_{L}}^{-1}\cdot X+ X^{\rm{T}}\cdot {B_{L}}^{-1} \cdot s\right)\ \ {\rm{mod}}\ 1.
\label{eq:surgery}
\end{equation}

\begin{example} \label{ex:example-algsplit}
Suppose that the surgery link $L$ is algebraically split
(ie, $B_L$ is diagonal). As before, denote by
$m_i$ the meridian of $L_i$ oriented so that
$\hbox{lk}_{\mathbf{S}^3}(L_i,m_i)=+1$ and let $[m_i]\in H$ be
its homology class in $M$. It follows from (\ref{eq:identification}) and
the orientation convention that
\begin{equation}
\label{eq:calculsplit}
\phi_{M,\sigma}([m_{i}]) =
-\frac{1}{2b_{ii}}(1 - s_{i})\ \hbox{mod}\ 1.
\end{equation}
\end{example}

\section{Proof of the Theorem} 

\label{sec:proof}

A technical difficulty lies in the computation of $q_{M,\sigma}$
from the torsion $\tau(M,\sigma)$. Fortunately, $\tau(M,\sigma)$
can be computed from a surgery presentation of $M$ and a charge
representing $\sigma$ (see \cite{TSurg} or \cite{TTorsions}). In the 
previous section, we computed $\phi_{M,\sigma}$ from a surgery
presentation of $M$ and a Chern vector representing $\sigma$.
Thus, the proof consists in two steps: 1. compare
charges to Chern vectors (there must be a bijective correspondence between them);
2. compare $q_{M,\sigma}$ to $\phi_{M,\sigma}$ using surgery presentations. \\

We shall use the notations of the previous section. In particular, we have fixed 
an ordered oriented framed $n$--component link $L$ in $\mathbf{S}^{3}$, 
such that the oriented $3$--manifold $V_L$ obtained
by surgery along $L$ is diffeomorphic to our oriented rational homology $3$--sphere $M$. \\

The comparison of the two combinatorial descriptions of Spin$^c(V_L)$
is contained in the following

\begin{claim}
\label{claim:comparaison}
If $\sigma \in {\rm{Spin}}^c(V_L)$ corresponds
to $[k] \in \mathcal{C}_{L}$,
then $\sigma$ corresponds to $[s]\in \mathcal{V}_{L}$, where
\begin{equation}
\label{eq:comparaison}
\forall j\in\{1,\dots,n\},\quad
s_j=1-k_j+\sum_{i=1}^n b_{ij}.
\end{equation}
\end{claim}

\begin{remark}
Claim \ref{claim:comparaison} is true for \emph{any} closed oriented
connected $3$--manifold (ie, with arbitrary first Betti number).
\end{remark}

\begin{proof}[Proof of the Claim \ref{claim:comparaison}]
We denote by $\sigma_2$ the distinguished relative \spinc--structure in
$\hbox{Spin}^c\left(\cup_{j=1}^{n} Z_j, {\cup}_{j=1}^{n} \sigma^0\right)$. Let also
$\sigma_1\in \hbox{Spin}^c({\mathbf{E}}, {\cup}_{j=1}^{n} \sigma^{0})$ be such that
$$\sigma=\sigma_1\cup \sigma_2\in \hbox{Spin}^c(V_L).$$ Pick an extension
$\tilde{\sigma}$ of $\sigma$ to $W_L$ and let
$\xi$ be the isomorphism class of $U(1)$--principal bundles determined by
$\tilde{\sigma}\in \hbox{Spin}^c(W_L)$.
On the one hand, the first Chern class $c_1(\xi)$ of $\xi$, when
expressed in the preferred basis $\left([S_j]^*\right)_{j=1}^n$
of $H^2\left(W_L\right)\simeq  \hbox{Hom}\left(H_2(W_L),\Z\right)$,
gives a multi--integer $s\in \Z^n$;
then $[s]\in \mathcal{V}_L$ corresponds
to $\sigma$. On the other hand, the Poincar\'e dual to the relative Chern class of
$\sigma_1\in \hbox{Spin}^c\left({\mathbf{E}}, {\cup}_{j=1}^n \sigma^{0}\right)$,
when expressed in the preferred basis $\left([m_j]\right)_{j=1}^n$
of $H_1({\mathbf{E}})$, gives a multi--integer $k\in \Z^n$;
then $[k]\in \mathcal{C}_{L}$ corresponds to $\sigma$.
Thus, proving that those specific integers $k$ and $s$ verify
(\ref{eq:comparaison}) modulo $2 \cdot\hbox{Im }B_L$ will be enough.

In the sequel we denote by $\left(\mathbf{S}^3\right)_{\varepsilon}$ a collar push-off
of $\mathbf{S}^3=\partial \mathbf{D}^4$ in the interior of  $\mathbf{D}^4$.
The surface $S_{j}$ can be decomposed (up to isotopy) in $W_L$ as
$$S_{j} = D_{j} \cup A_{j} \cup \left(-\Sigma_{j}^{\hbox{\scriptsize cut}}\right)_\varepsilon
\cup \bigcup_{l} \left(-R_{jl}\right)_\varepsilon$$
where the subsurfaces, illustrated on Figure \ref{fig:decomposing_Si}, 
are defined as follows:
\begin{figure}
\centerline{\relabelbox \small
\epsfxsize 4.4truein \epsfbox{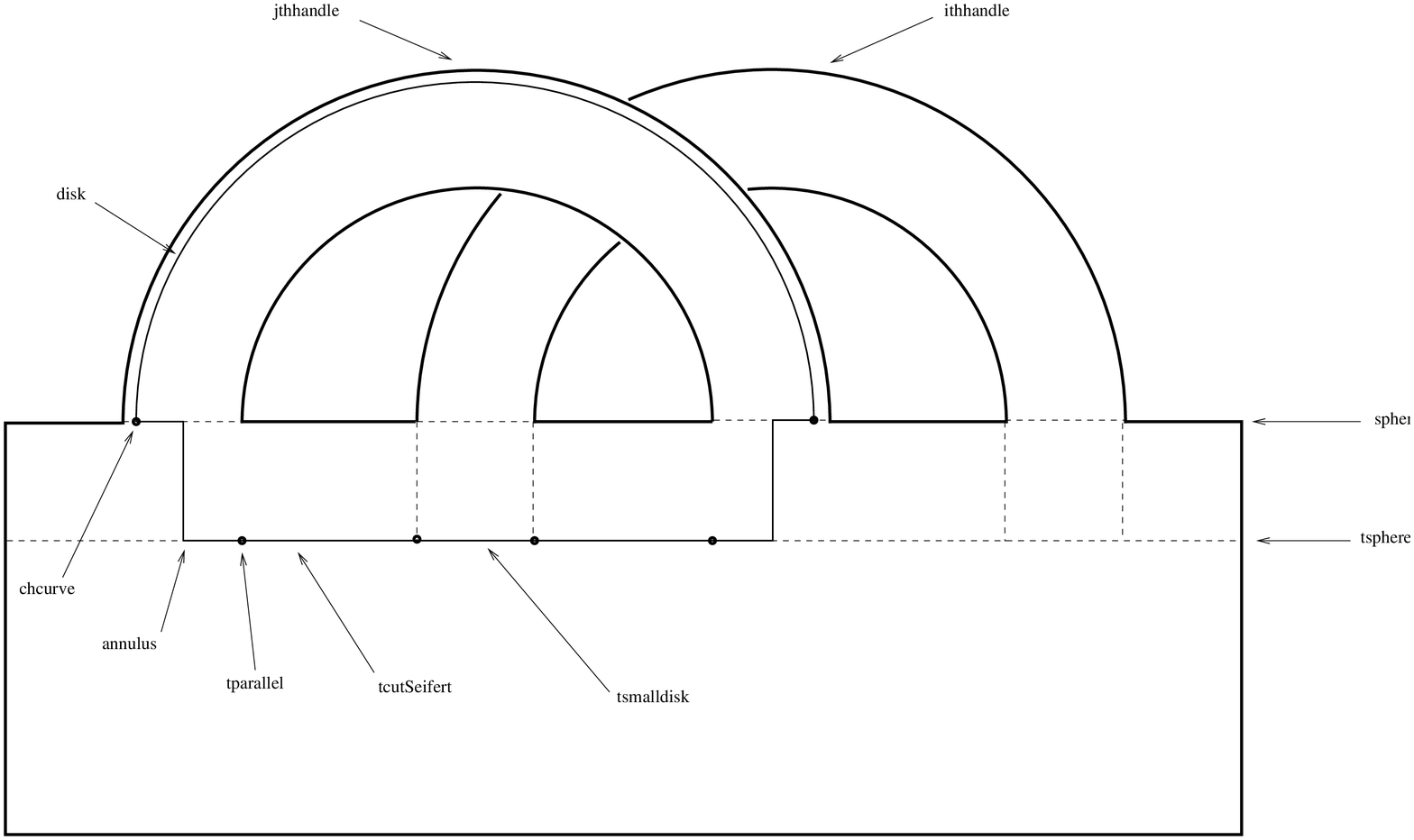}
\relabel  {ithhandle}{$i$--th handle}
\adjustrelabel <-1.2cm,0cm> {jthhandle}{$j$--th handle}
\adjustrelabel <-0.2cm,0cm> {disk}{$D_j$}
\adjustrelabel <0cm,-0.2cm> {chcurve}{$\lambda_j$}
\adjustrelabel <0.1cm,-0.2cm> {annulus}{$A_j$}
\adjustrelabel <0.1cm,-0.2cm> {tparallel}{$\left(l_j\right)_\varepsilon$}
\adjustrelabel <-0.1cm,-0.2cm> {tcutSeifert}{$\left(-\Sigma_j^{\hbox{\scriptsize cut}}\right)_\varepsilon$}
\adjustrelabel <-0.1cm,-0.2cm> {tsmalldisk}{one of the $\left(-R_{jl}\right)_\varepsilon$'s}
\adjustrelabel <0cm,-0.1cm> {tsphere}{$\left(\mathbf{S}^3\right)_\varepsilon$}
\adjustrelabel <0cm,-0.1cm> {sphere}{$\mathbf{S}^3$}
\endrelabelbox}
\caption{A decomposition of the surface $S_j$}
\label{fig:decomposing_Si}
\end{figure}
\begin{itemize}
\item[$\centerdot$] $D_{j}$ is a meridian disc of $Z_j$ such that
$\partial D_{j}$ is the characteristic curve
$\lambda_{j}$ of the $j$--th surgery;
\item[$\centerdot$] $A_{j}$ is the annulus of an isotopy of $-\lambda_{j}$ to $L_j$,
union the annulus of an isotopy of $-L_j$ to  $\left(L_j\right)_\varepsilon$,
union the annulus of an isotopy of $\left(-L_j\right)_\varepsilon$
to $\left(l_{j}\right)_\varepsilon$,
where $l_j$ denotes the preferred parallel of $L_j$ in $\mathbf{S}^3$
(ie, $\hbox{lk}_{\mathbf{S}^3}(l_j,L_j)=0$);
\item[$\centerdot$] $\Sigma_{j}$ is a Seifert surface for $l_j$ in ${\mathbf{S}}^{3}$
disjoint from $L_j$ and in transverse position with the $L_i$'s ($i\neq j$).
For each intersection point $x_{l}$ between $\Sigma_j$ and a $L_{i}$,
remove a small disc $R_{jl}$ so that
$\Sigma_j = \Sigma_{j}^{\hbox{\scriptsize cut}} \cup \bigcup_{l} R_{jl}$.
\end{itemize}
By definition of $s$, we have $s_{j} = \langle c_{1}(\xi), [S_{j}] \rangle =
\langle c_{1}(p|_{S_{j}}), [S_{j}] \rangle$ where $p$ is representative
for $\xi$ and where $c_{1}(p|_{S_{j}}) \in H^{2}(S_{j})$ is the obstruction
to trivialize $p$ over $S_{j}$.
So $P^{-1}c_{1}(p|_{S_{j}}) = s_{j}\cdot [\hbox{pt}] \in H_{0}(S_{j})$.
Let tr be a trivialization of $p$ on $\partial {\mathbf{E}}$ and let
$\hbox{tr}_\varepsilon$ be the corresponding
trivialization of $p$ on $\left(\partial {\mathbf{E}}\right)_\varepsilon$.
A classical argument (calculus of obstructions in compact oriented manifolds 
by means of  Poincar\'e dualities) leads to the equality
\begin{eqnarray}
\label{eq:sum_of_integers}
H_0(S_j)\ni P^{-1}c_{1}(p|_{S_{j}})& = &
i_{*}P^{-1}c_{1}\left(p|_{D_{j}}, \hbox{tr}|_{\lambda_{j}}\right)\\
&&+\ i_{*}P^{-1}c_{1}\left(p|_{A_{j}}, \hbox{tr}|_{-\lambda_{j}}
\cup \hbox{tr}_\varepsilon|_{(l_{j})_\varepsilon}\right) \nonumber\\
&&-\ i_{*}P^{-1}c_{1}\left(p|_{\left(\Sigma_{j}^{\hbox{\scriptsize cut}}\right)_\varepsilon},
\hbox{tr}_\varepsilon|_{\left(\partial \Sigma_{j}^{\hbox{\scriptsize cut}}\right)_\varepsilon}
\right) \nonumber\\
&&- \sum_{l} i_{*}P^{-1}c_{1}\left(p|_{\left(R_{jl}\right)_\varepsilon},
\hbox{tr}_\varepsilon|_{\left(\partial R_{jl}\right)_\varepsilon}\right),
\nonumber
\end{eqnarray}
where $P$ denotes a Poincar\'e duality isomorphism 
for the appropriate surface ($D_j$, $A_j$, $\Sigma_j^{\hbox{\scriptsize cut}}$ or $R_{jl}$). 
For an appropriate choice of $p$ in
the class $\xi$ and for an appropriate choice of tr, we have
\begin{eqnarray*}
c_1\left(p|_{{\mathbf{E}}},\hbox{tr}\right)&=&c(\sigma_1)\in H^2({\mathbf{E}},\partial {\mathbf{E}})\\
c_1\left(p|_{\cup_{j} Z_j},\hbox{tr}\right)&=&c(\sigma_2)
\in H^2\left(\cup_{j} Z_j,\cup_{j}\partial Z_j\right)\\
c_1\left(p|_{{\rm{N}}(L)},\hbox{tr}\right)&=&c(\sigma_3)
\in H^2\left({\rm{N}}(L),\partial {\rm{N}}(L)\right)
\end{eqnarray*}
where, in this last requirement, $\hbox{N}(L)$ is a tubular neighborhood of $L$
in $\mathbf{S}^3$ and $\sigma_3$ is an arbitrary element of
$\hbox{Spin}^c\left(\hbox{N}(L),\cup_{j} \sigma^{0}\right)$.
For such choices, we now compute separately
each term of the right hand side of (\ref{eq:sum_of_integers}).
\begin{enumerate}
\item The first term is of the form $d_{j}\cdot[\hbox{pt}]$. Here
$$d_{j} = \langle c(\sigma_2),[D_j]\rangle
= -\left(\textrm{oriented core of }Z_j\right) \bullet [D_{j}]=+1,$$ where
the intersection is taken in $Z_j$. (Note that $Z_j=\left(\mathbf{D}^2\times
\mathbf{S}^1\right)_j$ if we denote by $\left(\mathbf{D}^2\times
\mathbf{D}^2\right)_j$ the $2$--handle of $W_L$ corresponding to $L_j$,
and be careful of the fact that the above specified oriented core of $Z_j$
is $-\left(0\times \mathbf{S}^1\right)_j$.)
\item The second term is of the form $a_{j}\cdot[\hbox{pt}]$.
Here $a_{j} = \langle c(\sigma_3), [A_{j}] \rangle$ where
$A_{j}$ is regarded as a relative $2$--cycle in $(\hbox{N}(L), \partial \hbox{N}(L))$
once the collar has been squeezed. Since $\partial A_j$ is
$-\lambda_j \cup l_j$,
$[A_{j}]$ is $-b_{jj}$ times the class of the meridian disc of $L_j$
(oriented so that its oriented boundary is $m_j$)
in $H_2(\hbox{N}(L),\partial \hbox{N}(L))$. Then, $a_j=-b_{jj}\cdot \rho_{j}$ 
where $\rho_j$ is defined to be
$$\rho_{j}=\langle c(\sigma_3), [\hbox{meridian disc of $L_j$}]
\rangle\in\Z.$$ Note that $\rho_{j}\equiv J_{\sigma^{0}}([m_{j}]) \equiv 1$
mod $2$ (by the Atiyah-Johnson correspondence, see Lemma \ref{lem:Chern_mod2}).
\item The third term is $-g_{j}\cdot [\hbox{pt}]$ where 
$g_{j} = \langle c(\sigma_1),[\Sigma_{j}^{\hbox{\scriptsize cut}}] \rangle$. 
But, that integer is equal to
$$
g_j = \left(P^{-1}c(\sigma_1)\right) \bullet [\Sigma_j^{\hbox{\scriptsize cut}}] =
\left(\sum_{i} k_{i} [m_{i}] \right) \bullet [\Sigma_{j}^{\hbox{\scriptsize cut}}]=
\sum_{i} k_i\delta_{ij} = k_{j}
$$ 
where the intersection is taken in ${\mathbf{E}}$.
\item The fourth term is given by $-\sum_{l} r_{jl}\cdot [\hbox{pt}]$.
Here $r_{jl} = \langle c(\sigma_3), [R_{jl}] \rangle$.
For each index $l$, denote by $i(l)$ the integer $i$ such that $x_l$
is an intersection point of $\Sigma_j$ with $L_{i}$, and denote by
$\epsilon(l)$ the sign of the intersection point $x_l$.
Then, from the definition of $\rho_i$ (given for the second term),
we have $r_{jl}=\epsilon(l)\cdot\rho_{i(l)}$.
Hence
$$\sum_{l} r_{jl}=\sum_{\substack{i=1\\i\neq j}}^n
b_{ij} \rho_{i}.$$
\end{enumerate}
Putting those computations together, we obtain that (\ref{eq:sum_of_integers})
is equivalent to the identity
\begin{eqnarray*}
s_{j} &=& d_j + a_j -g_j - \sum_l r_{jl} \\
&=& 1 - b_{jj}\rho_j -k_j-\sum_{\substack{i=1\\i\neq j}}^n b_{ij}\rho_{i}\\
&=& \left(1 -k_j + \sum_{i=1}^n b_{ij}\right)
          -\sum_{i=1}^n b_{ij}\left(\rho_{i}+1\right).
\end{eqnarray*}
The claim now follows from the fact that $\rho_{i}\equiv 1$ mod $2$
for all $i=1,\dots,n$.
\end{proof}
We are now able to prove the Theorem.
Assume first that $M$ is obtained by surgery along
an algebraically split link $L$, and that $\sigma$ is represented
by a charge $k$ on $L$.
Then, according to \cite[Chapter X, \S 5.4]{TTorsions}, we have that
$$q_{M,\sigma}([m_{j}]) = \frac{1}{2} - \frac{k_{j}}{2b_{jj}}\ \hbox{mod 1}.$$
Substituting $k_j = 1 - s_{j} + \sum_{i} b_{ij}$, we find that
this formula agrees with (\ref{eq:calculsplit}) of Example
\ref{ex:example-algsplit}. This proves the Theorem in this
particular case. Now consider the general case, when $L$ is not
necessarily algebraically split. We shall use the following
observation due to Ohtsuki.
\begin{lemma} \label{algsplitsuffit}
Let $M$ be an oriented rational homology $3$--sphere. There exist
non-zero integers $n_{1}, \ldots, n_{r}$ such that $M \# L(n_{1},
1) \# \cdots \# L(n_{r}, 1)$ can be presented by surgery along a
framed link $L$ algebraically split in ${\mathbf{S}}^{3}$.
\end{lemma}
Here $\#$ denotes connected sum and $L(n,1)$ is the
$3$--dimensional lens space obtained by surgery along a trivial
knot with framing $n \not= 0$ in ${\mathbf{S}}^{3}$. Apply that lemma to
the oriented rational homology $3$--sphere $M$ we are working with, and consider 
the resulting manifold $M'= M \# L(n_{1}, 1) \# \cdots \# L(n_{r}, 1)$. 
Set $\sigma'=\sigma
\# \sigma_1 \# \cdots \#\sigma_r \in \hbox{Spin}^c(M')$ where
$\sigma_1,\dots, \sigma_r$ denote arbitrary \spinc--structures on
the lens spaces. Then, we have $q_{M',\sigma'} = \phi_{M',\sigma'}$.
By definition of $\#$, there is a small $3$--ball $B \subset M$ such
that $M \setminus B \subset M'$. This inclusion induces a
(injective) homomorphism $i_{*}:H_{1}(M) \to H_{1}(M')$. Since we
can compute $\phi_{M',\sigma'}$ from a split surgery presentation of $M'$ using
the surgery formula (\ref{eq:surgery}), we have that $\phi_{M,\sigma}
= \phi_{M',\sigma'} \circ i_{*}$. It follows from 
\cite[Chapter XII, \S 1.2]{TTorsions} (which describes the behaviour 
of the Reidemeister--Turaev torsion under $\#$) that, similarly,
$q_{M,\sigma}= q_{M',\sigma'} \circ i_{*}$. We deduce that 
$q_{M,\sigma}=\phi_{M,\sigma}$ and we are done.

\end{document}
